\documentclass[a4paper,12pt]{article}

\usepackage{latexsym}
\usepackage{amssymb,amsmath}

\usepackage{color}

\usepackage{scrtime}

\usepackage{hyperref}
\font\corsivo=rsfs10 at 12pt 
\font\doppio=msbm10 at 12pt 
\font\scdoppio=msbm8 

\newcommand{\R}{\hbox{\doppio R}}

\newcommand{\rn}{{\hbox{\scdoppio R}^N}}
\newcommand{\G}{\hbox{\doppio G}}     
\newcommand{\hei}{\hbox{\doppio H}} 
\newcommand{\RN}{{\R}^N}
\newcommand{\RR}{\hbox{\scdoppio R}}    

\newcommand{\A}{\mbox{\corsivo A}}      

\newcommand{\C}{\mbox{\corsivo C}}      

\newcommand{\grh}{\nabla_{\!\!H}}               
\newcommand{\lh}{\Delta_H}              

\newcommand{\Cuno}{\mbox{\corsivo C}^{\,1}}

\newcommand{\gr}{\nabla}
\newcommand{\grl}{\nabla_{\!\!L}}               
\newcommand{\diver}{\mathrm{div}}       
\newcommand{\diverl}{\diver_{\!L}}      
\newcommand{\decl}{:=}                  

\newcommand{\be}{\begin{equation}}      
\newcommand{\ee}{\end{equation}}        
\newcommand{\bern}{\begin{eqnarray*}}   
\newcommand{\eern}{\end{eqnarray*}}     

\newcommand{\sign}{\mathrm{sign}^+}

\newcommand{\norm}[1]{\abs{\abs{#1}}}

\newcommand{\SPC}{\mbox{{\bf S}-$p$-{\bf C}}}
\newcommand{\WPC}{\mbox{{\bf W}-$p$-{\bf C}}}

\newcommand{\MPC}{\mbox{{\bf M}-$p$-{\bf C}}}

\newtheorem{theorem}{Theorem}[section]
\newtheorem{lemma}[theorem]{Lemma}
\newtheorem{corollary}[theorem]{Corollary}
\newtheorem{definition}[theorem]{Definition}
\newtheorem{example}[theorem]{Example}
\newtheorem{remark}[theorem]{Remark}
\newtheorem{proposition}[theorem]{Proposition}

\newcommand{\bp}{\noindent{\bf Proof. }}
\newcommand{\ep}{\logend\medskip}
\newcommand{\logend}{\hspace*{\fill}$\Box$}
\newcommand{\abs}[1]{\left |#1\right |}
\newcommand{\mint}{{\int\!\!\!\!\!\!-}}

\newcommand{\plap}[1]{\diverl(\abs{\grl #1}^{p-2}\grl #1) }

\title{Uniqueness of $\sigma$-regular solutions \\ 
  of quasilinear elliptic problems}

\author{Lorenzo D'Ambrosio \\
{\small  Dipartimento di Matematica,
Universit\`a degli Studi di Bari}\\
{\small via E. Orabona, 4,  I-70125 Bari, Italy,
 {\tt dambros@dm.uniba.it}} \\ 
Enzo Mitidieri\\
{\small Dipartimento di Matematica e Informatica,
Universit\`a degli Studi di Trieste}\\
{\small via A.Valerio, 12/1, I-34127 Trieste, Italy,
 {\tt mitidier@units.it  }}
}
\date{{\small Rapporto del Dipartimento di Matematica di Bari 13/12, July 11,  2012.}}

\begin{document}

\maketitle

\begin{abstract}
We study the uniqueness problem of $\sigma$-regular solution of the equation, 
$$-\Delta_p u+ \abs u^{q-1}u =h  \quad on\quad  \RN, $$
where $q>p-1>0.$ and $N> p.$
Other  coercive type equations associated to more general
differential operators are also investigated. Our uniqueness results 
hold for equations associated to the mean curvature type operators
 as well as for more general quasilinear subelliptic operators.
\end{abstract}

\tableofcontents

\section{Introduction}

Nonlinear elliptic problems of coercive type is still  a subject of vital interest in the PDE  circles.
As it is well known, coercive problems have they roots in the  classical calculus of variations
and  precisely in the problems related to the  existence of minima for convex functional.

In a celebrated paper \cite{Boccardo}, Boccardo, Galluet and Vazquez  studied, among other things,  the simplest canonical
quasilinear problem with non regular data,
\begin{equation}\label{1}
   -\Delta_p u+ \abs u^{q-1}u =h  \quad on\quad  \RN, 
\end{equation}
where $q>p-1>0$ and  $h\in L^1_{loc}(\RN)$. 

An earlier and important contribution to this problem in the case $p=2$, was obtained by Brezis \cite{bre84}.  Indeed he  proved that  for any $h\in L^1_{loc}(\RN)$ the semilinear equation (\ref{1}) has a unique distributional solution $u\in D'(\RN)\cap L^q_{loc}(\RN).$

For the general case $p>1$ existence results have been obtained later in \cite{Boccardo}.  These Authors, by using a clever approximation procedure, proved that if $q>p-1$ and $\displaystyle{p> 2-\frac{1}{N}},$ then
for any $h\in L^1_{loc}(\RN)$ the equation (\ref{1}) possesses a solution belonging to the space $X=W^{1,1}_{loc}(\RN)\cap W^{1,p-1}_{loc}(\RN)\cap  L^q_{loc}(\RN)$.
No general results about uniqueness were claimed in that paper.

In this work,  we shall study the uniqueness problem of solutions of (\ref{1}) and related  qualitative properties. We emphasize that we shall prove the uniqueness of solutions
of (\ref{1}) in the space $W^{1,p}_{loc}(\RN)\cap L^q_{loc}(\RN).$
To this end, first we set up two essential tools which are of independent interest. Namely, the regularity of weak solutions of (\ref{1}) in the space $W^{1,p}_{loc}(\RN)\cap L^q_{loc}(\RN)$ and  comparison results on $\RN$ for related inequalities.
Further  we shall derive some property of the solutions of  (\ref{1}).

Our approach works also when dealing with more general operators and related inequalities. 
In this paper for sake of simplicity, we shall limit ourselves to  coercive problems in the  Carnot groups framework.
Clearly  this setting includes as special case  the Euclidean framework.

The approach we propose in this paper can be  successfully applied even when the differential  operator is not  the $p$-Laplacian operator.
Indeed the same uniqueness problem can be studied  for equations associated to  the mean curvature operator as well as  extensions of it.

The main results proved in this paper  are the following.

\begin{theorem} Let $1<p<2$, $q>1$, $h\in L^1_{loc}(\RN),$ then the problem 
$$ -\Delta_p u+ \abs u^{q-1}u =h  \quad on\quad  \RN,  $$
has at most one weak solution  $v\in W^{1,p}_{loc}(\RN)\cap L^q_{loc}(\RN)$.
 Moreover,   $$ \inf_{\RR^N} h\le \abs v^{q-1} v \le \sup_{\RR^N} h. $$
\end{theorem}
\begin{theorem} Let  $q>0$, $h\in L^1_{loc}(\RN)$ then the problem ,
$$ -\diver\frac{\gr u}{\sqrt{1+\abs{\gr u}^2}} + \abs u^{q-1}u =h  \quad on\quad  \RN,  $$
has at most one weak solution  $v\in W^{1,1}_{loc}(\RN)\cap L^q_{loc}(\RN)$.
 Moreover, 
$$ \inf_{\RR^N} h\le \abs v^{q-1} v \le \sup_{\RR^N} h. $$
\end{theorem}
With further assumptions on $h$ or on the solutions, we have the following results.
\begin{corollary}  Let $q>p-1>1$, let $h\in L^{\infty}(\RN)$. If 
  $$ 2<p<\frac{2N}{N-1},$$
  then the problem,
  $$ -\Delta_p u+ \abs u^{q-1}u =h  \quad on\quad  \RN,  $$
has at most one weak solution.
\end{corollary}
\begin{theorem}\label{14} Let $p>2$, $q>p-1.$  
  If  $h\in L^1_{loc}(\RN)$ then the problem,
  $$ -\Delta_p u+ \abs u^{q-1}u =h  \quad on\quad  \RN,  $$
has at most one weak solution $v$ in the class,
$$ \left\{ u\in W^{1,p}_{loc}(\RN)\cap L^q_{loc}(\RN) : there\ exist\ 
   \theta<\frac{1}{p-2}\ such\ that \ 
   \abs{\nabla u(x)} \le c \abs x^{\theta}\ for\ \abs x\  large \right\}. 
$$
 Moreover, 
$$ \inf_{\RR^N} h\le \abs v^{q-1} v \le \sup_{\RR^N} h. $$
\end{theorem}

Other partial results for the case $p>2$ are presented in Section \ref{sec:comp1}.

Our uniqueness results concern  solutions that belong to  the class $W^{1,p}_{loc}(\RN)\cap L^q_{loc}(\RN)$.  Of course, this set is contained in the space $X$ considered in \cite{Boccardo}. However we point out that,  when dealing with uniqueness results additional regularity  is usually required by several Authors.
See for instance \cite{Boccardo2}.
Indeed, in that paper the Authors obtain the existence of solutions of problem (\ref{1}) belonging to a certain space $T^{1,p}_0$.
The uniqueness result proved  in \cite{Boccardo2} concerns entropy solutions. 

We also emphasize that in this paper  we shall  restrict our analysis to the case $q>1.$ Indeed, as it is well known, see \cite{Boccardo2}, if  $1<p\le 2-(1/N)$ and $q=1, $ then there exists  $h\in L^1_{loc}(\RN)$ such that
(\ref{1}) has no solutions belonging to $W^{1,1}_{loc}(\RN)$. 
\medskip

The paper is organized as follow. In the next section we describe the setting and the notations.
In Section \ref{sec:est} we prove some a priori estimates on the solutions of the problems.
Section \ref{sec:comp} is devoted to prove the comparison results and to derive some of their consequences.
In this paper a preeminent role is played by the \MPC\ operators (see below for the definition).
In the appendix \ref{appMPC} we prove some inequalities that   guaranty that an operator is \MPC.
In Appendix \ref{appCarnot}, for the convenience of the reader,  we collect some basic facts about the Carnot groups.
\medskip

\noindent{\bf Note.} The results of this paper have been announced by the second author in 
Rome on May 7, 2012 and during a PDEs workshop dedicated to Patrizia Pucci's birthday in Perugia
  on May 30, 2012.
In the latter occasion,
as an outcome of several discussions with Professor James Serrin and Professor Alberto Farina,
we learned that they have obtained similar results to those proved in this paper.
An expanded  version of this work will appear in \cite{dfms}.

\section{Notations and  definitions}\label{sec:prehi}

In this paper $\nabla$ and $\abs\cdot$ stand  respectively for the usual gradient in $\RN$
and the Euclidean norm. $\Omega\subset \RN$ open.

Throughout this paper we shall use some concepts briefly described in the Appendix \ref{appCarnot}.
For further details related to Carnot groups the interested reader may refer to \cite{bon-lan-ugu07}.

\medskip

Let $\mu\in\C(\R^N;\R^l)$ be a matrix 
$\mu\decl(\mu_{ij})$, $i=1,\dots,l$, $j=1,\dots,N$
and assume that for any $i=1,\dots,l$, $j=1,\dots,N$ the derivative
 $\frac{\partial}{\partial x_j} \mu_{ij}\in\C(\Omega)$.
For $i=1,\dots,l$,  let $X_i$ and its formal adjoint $X_i^*$
be defined as
\be X_i\decl\sum_{j=1}^N \mu_{ij}(\xi)\frac{\partial}{\partial \xi_j},
\qquad X_i^*\decl-\sum_{j=1}^N \frac{\partial}{\partial\xi_j}\left(\mu_{ij}(\xi)\cdot\right),
   \label{mu}\ee
and let $\grl$ be the vector field defined by
$$\grl\decl (X_1,\dots,X_l)^T=\mu\nabla,$$
and $$\grl^*\decl(X_1^*,\dots,X_l^*)^T.$$

For any vector field $h=(h_1,\dots,h_l)^T\in\Cuno(\Omega,\R^l)$, we shall use
the following notation $ \diverl(h)\decl\diver\left(\mu^T h\right)$,
that is
\[ \diverl(h)=-\sum_{i=1}^l X_i^*h_i=-\grl^*\cdot h.\]

Examples of vector fields, which we are interested in, are the usual gradient
acting on $l(\le N)$ variables, vector fields
related to Bouendi-Grushin operator,
Heisenberg-Kohn sub--Laplacian,
Heisenberg-Greiner operator,
sub--Laplacian on Carnot Groups (see Appendix \ref{appCarnot}).

However in order to avoid cumbersome notations we shall limit ourselves to consider
the Carnot group case, which includes the Euclidean case.
We shall assume that the matrix $\mu$ has $\C^\infty$ entries
and that  $\RN$ can be endowed with a group law $\circ$ such that
the operator
 $$ L u\decl \diverl(\grl u),$$
is a sublaplacian on $(\RN, \circ)$. See Appendix \ref{appCarnot}
for further details.

Notice that 
in the Euclidean framework we have $\mu=I_N$, the identity matrix on $\RN$.

\bigskip

 In what follows we shall assume that
$\A:\RN\times\R\times\R^l\to\R^l$ is a Caratheodory
function, that is
for each $t\in\R$ and $\xi\in\R^l$ the function $\A(\cdot,t,\xi)$ is measurable; and
for a.e. $x\in\RN$, $\A(x,\cdot,\cdot)$ is continuous.

We consider operators $L$ ``generated'' by $\A$, that is
\[ L(u)(x)=\diverl \left(\A(x,u(x),\grl u(x))\right). \]
Our canonical model cases are the $p$-Laplacian operator, the mean curvature operator and
some related generalizations. See Examples \ref{examp:oper} below.

\begin{definition} \label{ptype} Let $\A:\RN\times\R\times\R^l\to\R^l$ be a 
  Caratheodory function.
  The function $\A$ is called \emph{weakly elliptic} if it generates a weakly elliptic operator $L$ i.e.
\[
\begin{array}{c}
  \A(x,t,\xi)\cdot \xi\ge 0\ \  \mathit{for\ each\ }x\in\RN,\,t\in\R,\,\xi\in\R^l,\\ \\
\A(x,0,\xi)=0\ \ \mathit{or}\ \  \A(x,t,0)=0
\end{array}
\eqno{(WE)}\]

Let $p\ge1$, the function $\A$ is called \WPC\ (weakly-$p$-coercive)
(see \cite{bid01}),
 if $\A$ is (WE) and it generates a weakly-p-coercive operator $L$, i.e.
 if there exists a constant
$k_2>0$ such that
\[ (\A(x,t,\xi)\cdot \xi)^{p-1}\ge k_2^{p} \abs{\A(x,t,\xi)}^{p}\ \mathrm{for\ each\ }x\in\RN,\,t\in\R,\,\xi\in\R^l.\eqno{(\WPC)}\]

Let $p>1$, the function $\A$ is called {\bf S-$p$-C} (strongly-$p$-coercive)
(see \cite{ser64, bid01, mit-poh01b}), 
if there exist $k_1,k_2 >0$ constants such that
\[ (\A(x,t,\xi)\cdot \xi)\ge k_1 \abs \xi^p\ge k_2^{p'} \abs{\A(x,t,\xi)}^{p'} \ \mathrm{for\ each\ }x\in\RN,\,t\in\R,\,\xi\in\R^l.\eqno{(\SPC)}\]
\end{definition}

\begin{definition} \label{def:sol} Let $\Omega\subset\RN$ be an open set and
 let $f:\Omega\times\R\times\R^l \to\R$ be a 
Caratheodory function.  Let $p\ge 1$.  We say that
  $u\in W^{1,p}_{loc}(\Omega)$ is a \emph{weak solution}  of
  $$ \diverl \left(\A(x,u,\grl u)\right)\ge f(x,u,\grl u)\qquad on\ \ \Omega,
   $$
 if   $\A(\cdot,u,\gr u)\in L^{p'}_{loc}(\Omega)$,
 $f(\cdot,u,\grl u)\in L^1_{loc}(\Omega)$, 
  and for any nonnegative $\phi\in\C^1_0(\Omega)$ we have
\[ - \int_{\Omega} \A(x,u,\grl u) \cdot\grl\phi \ge\int_{\Omega} f(x,u,\grl u)\phi. \]
\end{definition}

\begin{example}\label{examp:oper}
\begin{enumerate}
\item Let $p>1$. The  $p$-Laplacian  operator defined on suitable functions $u$ by,
$$\Delta_{p} u=\diverl\left(\abs{\grl u}^{p-2}\grl u\right)$$ 
is an operator generated by
$\A(x,t,\xi):=\abs \xi^{p-2}\xi$ which is  \SPC.

\item  If $\A$ is of mean curvature type, that is
$\A$ can be written as $\A(x,t,\xi):=A(\abs \xi) \xi$ with
  $A:\R\to \R$  a positive bounded continuous function
(see \cite{mit-poh01,bid01}),
then $\A$ is  {\bf W-$2$-C}.

\item The mean curvature operator in non parametric form 
$$Tu: =\diver\left(\frac{\nabla u}{\sqrt{1+\abs{\nabla u}^2}}\right),$$
is generated by $\A(x,t,\xi):=\frac{\xi}{\sqrt{1+\abs{\xi}^2}}$. In this case $\A$ is
  \WPC\ with $1\le p\le2$ and of mean curvature type but it is not  {\bf S-$2$-C}.
\item Let $m>1$. The operator
$$ T_mu\decl\diver\left(\frac{\abs{\gr u}^{m-2}\gr u}{\sqrt{1+\abs{\gr u}^m}}\right)
$$
is \WPC\ for $m\ge p\ge m/2$.

\end{enumerate}
\end{example}

\begin{definition} Let $\A:\RN\times \R^l\to \R^l$ be a Charateodory function.
  We say that $\A$ is monotone if 
  \be  (\A(x,\xi)-\A(x,\eta)) \cdot (\xi-\eta) \ge 0
   \qquad for\ \ \xi,\eta\in\R^l. \label{hypAmon}\ee
  Let $p\ge 1$. We say that $\A$ is \MPC\ (monotone $p$-coercive) if
  $\A$ is monotone and if there exists $k_2>0$ such that
  \begin{equation}
    \label{CWPC}
    \left( (\A(x,\xi)-\A(x,\eta))\cdot (\xi-\eta) \right)^{p-1}\ge k_2^{p}
    \abs{\A(x,\xi)-\A(x,\eta) }^{p}.
  \end{equation}
\end{definition}

\begin{example} 1. Let $1<p\le 2$ the function 
  $\A(\xi)\decl \abs {\xi}^{p-2}\xi$ is \MPC\ (see Appendix \ref{appMPC} for details).
  Therefore the following theorems apply to the $p$-Laplacian operator. 

  2. The mean curvature operator is \MPC\ with $1\le p\le 2$ (see appendix \ref{appMPC}).
\end{example}

\section{A priori estimates}\label{sec:est}

The following lemma  is a slight variation of a result proved in \cite{dam-mit:kato}. 
For easy reference  we shall include the detailed proof.

We shall consider  the following inequality,
\be \diverl\left(\A(x,v,\grl v)\right)-f\ge \diverl\left(\A(x,u,\grl u)\right)-g\qquad on\ \Omega.
\label{dis:comp0}\ee

\begin{theorem}\label{teo:start}  Let $\A:\Omega \times \R^l\to \R^l$ be \MPC,
  Let $f,g\in L^1_{loc}(\Omega)$ and let $(u,v)$ be weak solution of
  (\ref{dis:comp0}). Set $w\decl(v-u)^+$ and let $s>0$.
  If $(f-g) w\ge 0$ and $w^{s+p-1}\in L^1_{loc}(\Omega)$, then
  \be (f-g) w^s,\ \  (\A(x,\grl v)-\A(x,\grl u))\cdot \grl w\ w^{s-1}\in L^1_{loc}(\Omega)\label{regw+}.\ee
Moreover,  for any nonnegative $\phi\in\Cuno_0(\Omega)$ we have,
 \be  \int_\Omega (f-g) w^s \phi+c_1 s 
   \int_\Omega(\A(x,\grl v)-\A(x,\grl u))\cdot \grl w\ w^{s-1}\phi \le
     c_2 s^{1-p} \int_\Omega w^{s+p-1} \frac{\abs{\grl\phi}^p}{\phi^{p-1}},
      \label{dis:ests2i}\ee
  where $c_1= 1-\frac{p-1}{p}\left(\frac{\epsilon^{}}{k_2}\right)^{\frac{p}{p-1}}>0$, $c_2=\frac{p^p}{p\epsilon^p}$ and $\epsilon>0 $ is sufficiently small for $p>1$
 and $c_1=1$ and $c_2=1/k_2$ for $p=1$.
\end{theorem}

\begin{remark}\label{rem:wpc} i) Notice that from the above result it follows that 
  if  $u,v\in W^{1,p}_{loc}(\Omega)$ is a weak solution of
  (\ref{dis:comp0}), then $(f-g)\,w\in L^1_{loc}(\Omega)$.

 ii) The above lemma still holds if we replace the function 
  $f-g\in  L^1_{loc}(\Omega)$ with a regular Borel measure on $\Omega$.

 iii)  If $(u,v)$ is a weak solution of  (\ref{dis:comp0}) 
  and $u$ is a constant i.e.  $u\equiv const$, then  Theorem \ref{teo:start}
  still holds even for \WPC\ operators. See the following lemma.
\end{remark}

\begin{lemma}\label{lem:1sti} Let $\A$ be \WPC.
  Let $f,g\in L^1_{loc}(\Omega)$   and 
  let $v\in W^{1,p}_{loc}(\Omega)$ be a weak solution of
  \be \diverl\left(\A(x,u,\grl u)\right)\ge f-g,\qquad \ on\ \ \Omega. \label{dis:genloc}\ee
  Let $k>0$ and set $w\decl(v-k)^+$ and let $s>0$.
  If $(f-g) w\ge 0$ and $w^{s+p-1}\in L^1_{loc}(\Omega)$, then
  \be (f-g) w^s,\ \  \A(x,v,\grl v)\cdot \grl w\ w^{s-1}\in L^1_{loc}(\Omega)\label{reg}\ee
 and for any nonnegative $\phi\in\Cuno_0(\Omega)$ we have,
 \be  \int_\Omega (f-g) w^s \phi+c_1 s 
   \int_\Omega \A(x,v,\grl v)\cdot \grl w\ w^{s-1}\phi \le
     c_2 s^{1-p} \int_\Omega w^{s+p-1} \frac{\abs{\grl\phi}^p}{\phi^{p-1}},
      \label{dis:ests}\ee
  where $c_1$ and $c_2$ are as in Theorem \ref{teo:start}.
\end{lemma}

The above result lies on the following result proved in 
 \cite[Theorem 2.7]{dam-mit:kato}.
\begin{theorem}[\cite{dam-mit:kato}]\label{th:katocompare} 
 Let $\A:\Omega \times\R\times \R^N\to \R^N$ be monotone Caratheodory function. 
Let $f,g\in L^1_{loc}(\Omega)$ and let $u,v$ be weak solution of
\be \diverl\left(\A(x,v,\grl v)\right)-f\ge \diverl\left(\A(x,u,\grl u)\right)-g\qquad on\ \Omega.
\label{dis:comp3}\ee

Let $\gamma\in\Cuno(\R)$ be 
such that $0\le \gamma(t),\gamma'(t)\le M$, then
\begin{eqnarray}
 \lefteqn{- \int_\Omega(\A(x,v,\grl v)-\A(x,u,\grl u)) \cdot \grl \phi\  \gamma(v-u) \ge} \\
&&\ge \int_\Omega \gamma'(v-u)\ (\grl v-\grl u )\cdot
    \left(\A(x,v,\grl v)-\A(x,u,\grl u\right)\phi
\\&&\qquad +\int_\Omega \phi\gamma(v-u) (f-g)\ \ \ on\ \Omega.
\end{eqnarray}
Hence
$$ \diverl\left(\gamma(v-u)(\A(x,v,\grl v)-\A(x,u,\grl u)) \right)
    \ge \gamma(v-u) (f-g)\ \ \ on\ \Omega.$$
Moreover\footnote{We recall that the function $\sign$ is defined as  
   $\sign(t)\decl 0$ if $t\le 0$ and $\sign(t)\decl 1$ otherwise.}
\be \diverl\left(\sign(v-u)(\A(x,v,\grl v)-\A(x,u,\grl u)) \right)
    \ge \sign(v-u) (f-g)\ \ \ on\ \Omega.\label{dis:katocomp}
\ee
\end{theorem}
\noindent{\bf Proof of Theorem \ref{teo:start}.} 
  Let $\gamma\in\Cuno(\R)$ be a bounded nonnegative function with
  bounded nonnegative first derivative and let $\phi\in\Cuno_0(\Omega)$
  be a nonnegative test function.

  For simplicity we shall omit the arguments of $\A$. So we shall write
  $\A_u$  and $\A_v$ instead of $\A(x,  \grl u)$ and $\A(x,  \grl v)$
  respectively.

  Applying Lemma \ref{th:katocompare}, we obtain
  \begin{eqnarray}
  \int_\Omega (f-g)\gamma(w) \phi+ \int_\Omega(\A_v-\A_u)\cdot\ \grl w\ \gamma'(w)\phi&\le& - \int_\Omega{(\A_v-\A_u)}\cdot{\grl\phi}\ \gamma(w)\nonumber\\
    &\le& \int_\Omega\abs{\A_v-\A_u }\,\,\abs{\grl\phi} \gamma(w)\label{tec1}
  \end{eqnarray}

 Let $p>1$. From (\ref{tec1}) we have
\begin{eqnarray*}\lefteqn{\int_\Omega (f-g)\gamma(w) \phi+ \int_\Omega(\A_v-\A_u)\cdot\ \grl w\ \gamma'(w)\phi\le}\\
   &&\le\left(\int_\Omega\abs{\A_v-\A_u}^{p'}\gamma'(w)\phi\right)^{1/p'}
    \left(\int_\Omega\frac{\gamma(w)^p}{\gamma'(w)^{p-1}}  \frac{\abs{\grl\phi}^p}{\phi^{p-1}} \right)^{1/p}\\
 &&\le \frac{\epsilon^{p'}}{p'k_2^{p'}} \int_\Omega(\A_v-\A_u)\cdot \grl w\ \gamma'(w)\phi + 
   \frac{1}{p\epsilon^p}\int_\Omega\frac{\gamma(w)^p}{\gamma'(w)^{p-1}} \frac{\abs{\grl\phi}^p}{\phi^{p-1}},
  \end{eqnarray*}
  where $\epsilon>0$ and all integrals are well defined provided
  $ \frac{\gamma(w)^p}{\gamma'(w)^{p-1}} \in L^1_{loc}(\Omega)$. 
  With a suitable choice of $\epsilon>0,$ for any 
  nonnegative $\phi\in\Cuno_0(\Omega)$  and 
   $\gamma\in\Cuno(\R)$
  as above such that $\frac{\gamma(w)^p}{\gamma'(w)^{p-1}} \in L^1_{loc}(\Omega),$ 
 it follows that,
  \be  \int_\Omega (f-g)\gamma(w) \phi+c_1 \int_\Omega(\A_v-\A_u)\cdot \grl w\ \gamma'(w)\phi \le
     \frac{1}{p\epsilon^p} \int_\Omega\frac{\gamma(w)^p}{\gamma'(w)^{p-1}}  \frac{\abs{\grl\phi}^p}{\phi^{p-1}}.
      \label{primastima}\ee

  Now for $s> 0$, $1>\delta>0$ and $n\ge 1$,  define
  \be \gamma_n(t)\decl
  \begin{cases}
    (t+\delta)^s&if\ 0\le t<n-\delta,\\ \\
    \displaystyle cn^s-\frac{s}{\beta-1}{n^{\beta+s-1}}{(t+\delta)^{1-\beta}}&if\ t\ge n-\delta,\\
  \end{cases}  \label{gn}\ee
  where  $c\decl\frac{\beta-1+s}{\beta-1}$ and $\beta>1$ will be chosen later.
  Clearly $\gamma_n\in \C^1$,
  $$\gamma'_n(t)=
  \begin{cases}
    s(t+\delta)^{s-1}&if\ 0\le t<n-\delta,\\ \\
    s{n^{\beta+s-1}}{(t+\delta)^{-\beta}}&if\ t\ge n-\delta,\\
  \end{cases}  $$
  and  $\gamma_n$, $\gamma'_n$ are nonnegative and bounded  with 
  $\norm{\gamma_n}_\infty=cn^s$ and $\norm{\gamma_n'}_\infty=sn^{s-1}$. Moreover 
$$\frac{\gamma_n(t)^p}{\gamma_n'(t)^{p-1}}=
\begin{cases}
  s^{1-p}(t+\delta)^{s+p-1}&\ \  for\ t<n-\delta,\  \\ \\
  \theta(t,n) &\ \ for \ t\ge n-\delta,
\end{cases}
$$
where 
$$\theta(t,n)\decl\frac{(cn^s-\frac{s}{\beta-1}{n^{\beta+s-1}}{(t+\delta)^{1-\beta}})^p}{(s{n^{\beta+s-1}}{(t+\delta)^{-\beta}})^{p-1}}
 \le {(cn^s)^p}{s^{1-p}} n^{-(\beta+s-1)(p-1)}\, (t+\delta)^{\beta(p-1)}.$$
Choosing $\beta\decl \frac{s+p-1}{p-1}$ we have $c=p$, and 
$$\theta(t,n)\le p^p s^{1-p} n^{s p-(\beta+s-1)(p-1)}  (t+\delta)^{s+p-1}= p^p s^{1-p} (t+\delta)^{s+p-1}.$$
Therefore,  for $t\ge 0$ we have,
$$\frac{\gamma_n(t)^p}{\gamma_n'(t)^{p-1}}\le p^p s^{1-p} (t+\delta)^{s+p-1}.$$

Since by assumption  $w^{s+p-1}\in L^1_{loc}(\Omega)$, from (\ref{primastima})
 with $\gamma=\gamma_n$, it follows that
  $$\int_\Omega (f-g)\gamma_n(w) \phi+c_1 \int_\Omega(\A_v-\A_u)\cdot \grl w\ \gamma_n'(w)\phi \le
     \frac{p^p s^{1-p}}{p\epsilon^p}  \int_\Omega (w+\delta)^{s+p-1}\frac{\abs{\grl\phi}^p}{\phi^{p-1}}.$$

 Now, noticing that  $\gamma_n(t)\to (t+\delta)^s$ and $\gamma_n'(t)\to s (t+\delta)^{s-1}$ as 
  $n\to +\infty,$   $(f-g)\gamma_n\ge 0$ and $\A$ is monotone 
  (that is $(\A_v-\A_u)\cdot\grl w\ge 0$),  by Beppo Levi theorem 
  we obtain
  $$  \int_\Omega (f-g)\, (w+\delta)^s \phi+c_1 s \int_\Omega(\A_v-\A_u)\cdot \grl w\ (w+\delta)^{s-1}\phi \le
      c_2 s^{1-p}\int_\Omega (w+\delta)^{s+p-1} \frac{\abs{\grl\phi}^p}{\phi^{p-1}}.$$
By letting $\delta\to 0$ in the above inequality, we
 complete the proof of the claim in the case $p>1$.

Let $p=1$. From (\ref{tec1}) and the fact that $\A_v-\A_u$ is bounded,
 the estimate (\ref{primastima}) holds provided
we replace $p$ with $1$ and $\epsilon$ with $k_2$.
The remaining argument  is similar to the case $p>1$ and we omit it. 
\ep

\begin{remark}\label{rem:smin1}
  The assumption $u^{s+p-1}\in L^1_{loc}(\Omega)$, 
  is not needed for the statement (\ref{reg}). 
  Indeed what  really matters for the validity of (\ref{reg}) is the assumption  $u^{s+p-1}\in L^1_{loc}(S).$ Here $S$ is the support of $\grl \phi$.
  This remark will be useful when dealing with inequalities on unbounded set.
\end{remark}
\begin{lemma}\label{lem:sti2} Let $p\ge 1$ and let  $\A:\Omega \times \R^l\to \R^l$ be \MPC.
  Let $f,g\in L^1_{loc}(\Omega)$ and let $(u,v)$ be weak solution of
  (\ref{dis:comp0}). Set $w\decl(v-u)^+$.
  If $(f-g) w\ge 0$ and $w^{q}\in L^1_{loc}(\Omega)$ for $q>p-1$, then
  \be (f-g) w^{q-p+1},\ \  ((\A(x,\grl v)- \A(x,\grl u))\cdot \grl w\ w^{q-p}\in L^1_{loc}(\Omega),\label{reg+}\ee
 and for any $\varphi\in\Cuno_0(\Omega)$ such that $0\le\varphi\le 1$, we have,
\be\int (f-g)\sign(w) \ \varphi^\sigma\le c_3 
  \left(\frac{1}{\abs S}\int_S w^{q}\varphi^\sigma\ \right)^{\frac{p-1}{q}} 
    \left(\frac{1}{\abs S}\int_S \abs{\grl \varphi}^\sigma \right)^{\frac{p}{\sigma}} 
    \abs S,\label{glu}\ee
  where $S$ is the support of $\grl \varphi$,
  $c_3\decl \frac{\sigma^p}{k_2s^{p-1}}\left(\frac{c_2}{c_1}\right)^{(p-1)/p}$
  with
  $\sigma\ge \frac{pq}{q-p+1-s}$, $0<s<\min\{1,q-p+1\}$
  and $c_1, c_2$ as in the above Theorem \ref{teo:start}.
\end{lemma}
\bp The claim (\ref{reg+}) follows directly applying Theorem \ref{teo:start}.

 Let $s>0$ be such that $q\ge s+p-1$. 
  From Lemma \ref{teo:start} for any nonnegative $\phi\in\Cuno_0(\Omega)$
  we have  
   \be  \int (f-g) w^s \phi+c_1 s \int(\A_v -\A_u)\cdot \grl w\ w^{s-1}\phi \le
     c_2 s^{1-p} \int w^{s+p-1} \frac{\abs{\grl\phi}^{p}}{\phi^{p-1}},
     \label{est2}\ee
  where, as in the proof of Theorem \ref{teo:start}, we write 
  $\A_v$ and $\A_v$ for $\A(x,\grl v)$ and $\A(x,\grl u)$ respectively.

  Next, an application of Theorem \ref{th:katocompare} gives
(\ref{dis:katocomp}). That is 
\be \diverl\left(\sign(v-u)(\A(x,v,\grl v)-\A(x,u,\grl u)) \right)
    \ge \sign(v-u) (f-g)\ \ \ on\ \Omega.
\ee

  Next consider the case $p>1$. Let $0<s <\min\{ 1, q-p+1\}$.
  By definition of weak solution and 
   H\"older's inequality with exponent $p'$,
  taking into account that $\A$ is \MPC\ and from  (\ref{est2}) we get,
  \begin{eqnarray} \lefteqn{
    \int \sign w (f-g)\ \phi \le \int_S \abs{\A_v-\A_u}\abs{\grl \phi}\sign w}\label{dis:tec2}\\
     &&=
      \int_S\abs{\A_v-\A_u}w^{\frac{s-1}{p'}}\phi^{\frac{1}{p'}}\  
      \abs{\grl \phi}w^{\frac{1-s}{p'}}\phi^{-\frac{1}{p'}}
\\
  &&\le \frac{1}{k_2^{}} \left(\int_S  (\A_v-\A_u)\cdot\grl w\, w^{s-1} \phi\right)^{1/p'}
    \left(\int_S w^{(1-s)(p-1)}\frac{\abs{\grl\phi}^{p}}{\phi^{p-1}} \right)^{1/p}\\
  &&\le  \frac{1}{k_2^{}}\left(\frac{c_2}{c_1s^p}\right)^{1/p'}\left(\int_S w^{s+p-1} \frac{\abs{\grl\phi}^{p}}{\phi^{p-1}}\right)^{1/p'} \left(\int_S w^{(1-s)(p-1)}\frac{\abs{\grl\phi}^{p}}{\phi^{p-1}} \right)^{1/p}\label{est4}
     \end{eqnarray}
  Since  $q>s+p-1$ and $q>p-1$,  applying H\"older
  inequality to (\ref{est4}) with exponents
  $\chi\decl \frac{q}{s+p-1}$ and 
  $y\decl\frac{q}{(1-s)(p-1)}$, we obtain 
  \be
    \int \sign w \  (f-g)\phi \le c_3' 
\left(\int_S w^{q}\phi \right)^{\delta} 
    \left(\int   \frac{\abs{\grl\phi}^{p\chi'}}{\phi^{p\chi'-1}}\right)^{\frac{1}{p'\chi'}} 
\left(\int   \frac{\abs{\grl\phi}^{py'}}{\phi^{py'-1}}\right)^{\frac{1}{py'}}, \label{est5pre}\ee 
  where 
  $$\delta\decl \frac{1}{\chi p'}+\frac{1}{yp}=\frac{p-1}{q},\qquad
   c_3'\decl \left(\frac{c_2}{c_1s^p}\right)^{1/p'}\frac{1}{k_2}.$$
  Next for $\sigma\ge p\chi'$ (and hence 
  $\sigma> py'$  since $ p\chi'> py'$ )
 we choose $\phi\decl \varphi^\sigma$
  with $\varphi\in\Cuno_0(\Omega)$  such that $0\le\varphi\le 1$.
  Setting $S\decl support(\varphi)$, from (\ref{est5pre}) it follows that
    \be \int \sign w\  (f-g) \varphi^\sigma \le c_3' \sigma^p 
   \left(\int w^{q}\varphi^\sigma \right)^{\delta} 
   \left(\frac{1}{\abs S}\int_S  \abs{\grl\varphi}^{\sigma} 
        \right)^{\frac{p}{\sigma}}  \abs S^{1-\delta},\label{est5}\ee
completing the proof of (\ref{glu}).

Now, we assume that $p=1$. From (\ref{dis:tec2}), with the choice
  $\phi\decl \varphi^\sigma$, with $\varphi\in\Cuno_0(\Omega)$  such that $0\le\varphi\le 1$ and $\sigma \ge 1$, we have
  $$\int \sign w\ (f-g)\varphi^\sigma \le \frac{\sigma}{k_2}\int_S \abs{\grl\varphi}
  \le \frac{\sigma}{k_2}\left(\frac{1}{\abs S} \int_S \abs{\grl\varphi}^\sigma \right)^{1/\sigma},$$
  which concludes the proof.
\ep

Now, by specializing $f$ and $g$, we study
\be \diverl\left(\A(x,\grl v)\right)-\abs v^{q-1} v\ge \diverl\left(\A(x,\grl u)\right)-\abs u^{q-1} u\qquad on\ \Omega.
\label{dis:compq}\ee

\begin{lemma}\label{lem:estR}
  \begin{enumerate}
  \item Let $p>1$. Let   $\A$ be \MPC\  and
   let $q>\max\{1,p-1\}$.
  For any $\sigma>0$ large enough, there exists a constant 
  $c=c(\sigma,q,p,\A)>0$ such that if $(u,v)$ is weak solution of (\ref{dis:compq})
then for any nonnegative $\varphi\in \Cuno_0(\Omega)$ such that 
$\norm \varphi_\infty\le 1$, we have
\begin{eqnarray}
 \int (\abs v^{q-1}v-\abs u^{q-1}u)^+\varphi^\sigma\le c 
  \left(\frac{1}{\abs S}\int_S (v-u)^{+q}\varphi^\sigma\ \right)^{\frac{p-1}{q}} 
    \left(\frac{1}{\abs S}\int_S \abs{\grl \varphi}^\sigma \right)^{\frac{p}{\sigma}} 
    \abs S\label{vlu}\\
\int (\abs v^{q-1}v-\abs u^{q-1}u)^{+}\varphi^\sigma\le c\ \abs S\ 
   \left(\frac{1}{\abs S}\int_S \abs{\grl \varphi}^\sigma \right)^{\frac{1}{\sigma}\,\frac{p q}{q-p+1}}  \label{estv}
\end{eqnarray}
where $S\decl support( \varphi)$.

In particular if $B_{2R}\subset\subset\Omega$, then
\begin{equation}
  \label{dis:v-u}
  \left(\mint_{B_R} (\abs v^{q-1}v-\abs u^{q-1}u)^{+}\right)^{1/q} \le c\  R^{-\frac{p }{q-p+1}}.
\end{equation}

Moreover, for $x\in \Omega$, set $R=dist(x,\partial\Omega)/2$, we have
\begin{equation} \label{dis:v-ux}
\left(\mint_{B_R(x)} (\abs v^{q-1}v-\abs u^{q-1}u)^{+} \right)^{1/q}\le c\  dist(x,\partial\Omega)^{-\frac{p }{q-p+1}}. \end{equation}

\item Let $p=1$. Let   $\A$ be \MPC\  and
   let $q>0$.
  For any $\sigma>0$ large enough, there exists a constant 
  $c=c(\sigma,q,p,\A)>0$ such that if $(u,v)$ is weak solution of (\ref{dis:compq})
then for any nonnegative $\varphi\in \Cuno_0(\Omega)$ such that 
$\norm \varphi_\infty\le 1$, we have
\begin{eqnarray}
 \int (\abs v^{q-1}v-\abs u^{q-1}u)^+\varphi^\sigma\le c 
  \left(\frac{1}{\abs S}\int_S (v-u)^{+q}\varphi^\sigma\ \right)^{\frac{p-1}{q}} 
    \left(\frac{1}{\abs S}\int_S \abs{\grl \varphi}^\sigma \right)^{\frac{p}{\sigma}} 
    \abs S\label{vlu1}
\end{eqnarray}
where $S\decl support( \varphi)$.

In particular if $B_{2R}\subset\subset\Omega$, then
$$\left(\mint_{B_R} (\abs v^{q-1}v-\abs u^{q-1}u)^{+}\right)^{1/q} \le c\  R^{-\frac{1 }{q}}.
$$

Moreover,
 for $x\in \Omega$, set $R=dist(x,\partial\Omega)/2$, we have
 \begin{equation}
\left(\mint_{B_R(x)} (\abs v^{q-1}v-\abs u^{q-1}u)^{+} \right)^{1/q}\le c\  dist(x,\partial\Omega)^{-\frac{1}{q}}.
 \end{equation}
  \end{enumerate}
\end{lemma}

\bp Let $p>1$.
  From Lemma \ref{lem:sti2} we immediately obtain (\ref{vlu}).

  Reminding the well known inequality
  \be t^q-s^q\ge c_q (t-s)^q, \quad  for\ \  t>s\qquad(q>1),\label{wki}\ee
  from  (\ref{vlu}), we get  (\ref{estv}).

  In order to obtain the estimate (\ref{dis:v-u}) we specialize the test 
  function $\varphi$. Indeed, let $\phi\in\Cuno_0(\R)$ be such that
  $0\le \phi\le 1$, $\phi(t)=0$ if $\abs t\ge 2$ and 
  $\phi(t)=1$ if $\abs t\le 1$. Next,  we define $\phi_R(t)\decl\phi(t/R)$.
  The claim follows by choosing $\varphi(x)=\phi_R(\abs {x})$. 
  Indeed, with this choice we have $\abs{\grl \varphi}\le c R^{-1}$
  and $\abs S=\abs {A_R}=c R^N$.
  
  The estimate (\ref{dis:v-ux}) follows by choosing
  $\varphi(y)=\phi_R(\abs {y-x})$.

  The case $p=1$ follows the same argument as above so that  we can leave the details to the interested reader.
\ep

\begin{lemma}\label{lem:boot}  
  Assume that  either one of the following holds
  \begin{enumerate}
  \item Let $p> 1$, let  $\A$ be \MPC and $q>\max\{1,p-1\}$.
  \item Let $p= 1$, let  $\A$ be \MPC and $q> 0$.
  \end{enumerate}

 Let  $(u,v)$ be weak solution of (\ref{dis:compq}).
   Then $(v-u)^+\in L^r_{loc}(\Omega)$ for any $r<+\infty$.
\end{lemma}
\bp  Let $(u,v)$ be a solution of (\ref{dis:compq}) and set  $w\decl (v-u)^+$.
  Now since $w\in L^q_{loc}(\Omega)$, and  inequality \ref{wki} holds we are
  in the position to apply Theorem \ref{teo:start}, with $s=q-p+1$ obtaining
  $w^{q_1}\in L^1_{loc}(\Omega)$ with $q_1\decl{2q-p+1}$. 
  Applying again Theorem \ref{teo:start}, with $s=q_1-p+1$, we get
  $w^{q_2}\in L^1_{loc}(\Omega)$ with $q_2\decl q_1+q-p+1=q+2(q-p+1)$.
  Iterating $j$ times we have that
  $w^{q_j} \in L^1_{loc}(\Omega)$ with $q_j\decl q+j(q-p+1)$.
  Letting $j\to +\infty$ we have the claim.
\ep

\section{Comparison and Uniqueness}\label{sec:comp}

\begin{theorem}\label{teo:comp}  Let $p\ge 1$ and let $\A$ be \MPC.
  Let (u,v) be a $\sigma$-regular solution of
  \be \diverl\left(\A(x,\grl v)\right)-\abs {v}^{q-1}v
    \ge \diverl\left(\A(x,\grl u)\right)-\abs{u}^{q-1}u\qquad on\ \RN.
    \label{dis:cmp}\ee
  Assume that one of the following holds
  \begin{enumerate}
   \item Let $p> 1$,  and $q>\max\{1,p-1\}$.
   \item Let $p= 1$,  and $q> 0$.
  \end{enumerate} 
  Then $v\le u$ a.e. on $\RN$.
\end{theorem}

\bp Let $(u,v)$ be a solution of (\ref{dis:cmp}) and set  $w\decl (v-u)^+$.
  From Lemma \ref{lem:boot} we know that  $w\in L^r_{loc}(\RN)$ for any $r$,
  and hence we are
  in the position to apply Theorem \ref{teo:start} with $s$ large enough.
  Thus, from (\ref{wki}) and (\ref{dis:ests2i}) we get
$w^{q+s}\in L^1_{loc}$
$$\int w^{q+s} \phi \le c(s,q,p) \int w^{s+p-1} \frac{\abs{\grl\phi}^p}{\phi^{p-1}}.
$$
Applying the H\"older inequality with exponent $x\decl\frac{q+s}{s+p-1}>1$
we have
$$\int w^{q+s} \phi \le c(s,q,p) \int \frac{\abs{\grl\phi}^{px'}}{\phi^{px'-1}}.
$$
By the same choice  of  $\phi$ we made in Lemma \ref{lem:estR}, we have that
$$\int_{B_R} w^{q+s}  \le c R^{N-px'}=c R^{N-p(q+s)/(q-p+1)}.$$
Choosing $s$ large enough and letting $R\to +\infty$, we have that $w\equiv 0$ a.e.
that is the claim.
\ep

\begin{corollary}\label{cor:bd} Let $p\ge1$,
   let $\A$ be \WPC\ such that $\A(x,0)=0$. 
   Let $q>0$ be  as in Theorem \ref{teo:comp}.
Let  $h\in L^1_{loc}(\RN)$. 
   Let $v$ be a solution of the problem 
   \be -\diverl\left(\A(x,\grl v)\right)+\abs {v}^{q-1}v= h.\label{eq}\ee
  Then,
    $$ \inf_{\RR^N} h  \le\abs{ v}^{q-1} v\le \sup_{\RR^N} h.$$
  In particular, if $h\ge 0$ [resp. $\le 0$], then $v\ge 0$ [resp. $\le 0$] and
 if $h\in L^\infty(\RN)$, then $v \in L^\infty(\RN)$.

\end{corollary}
\bp  We prove one of the estimates, the remain one is similar.
  If $\sup_{\RR^N} h =+\infty$ there is nothing to prove. 
  Let $M\decl\sup_{\RR^N} h <+\infty$. 
  We define $u\decl \mathrm{sign}(M) \abs{M}^{1/q} $.
  Then
  $$ \diverl\left(\A(x,\grl v)\right)-\abs {v}^{q-1}v+h=0
  \ge h-M=
  \diverl\left(\A(x,\grl u)\right)-\abs{u}^{q-1}u-h, $$
  that is $(u,v)$ satisfy (\ref{dis:cmp}) with $u$ constant. 
   In this case  all the previous estimates still hold
   since in this case the operator can be seen as it were \MPC.
   See also Remark \ref{rem:wpc} and Lemma \ref{lem:1sti}.

 Thus the claim follows from 
  Theorem \ref{teo:comp}.
\ep
\begin{corollary}\label{cor:uniq} Let $p\ge1$ and
   let $\A$ be \WPC.
   Let $q>0$ be  as in Theorem \ref{teo:comp}.
Let  $h\in L^1_{loc}(\RN)$. 
   Then the eventual solution of  the problem (\ref{eq}) is unique.

   Moreover  if $\A(x,0)=0$ 
    and $v$ is a solution of (\ref{eq}), then  
 $$ \inf_{\RR^N} h  \le\abs{ v}^{q-1} v\le \sup_{\RR^N} h.$$
\end{corollary}
\bp Uniqueness. Let $u$ and $v$ two solutions of (\ref{eq}).
  Then $(u,v)$ solves 
  $$  \diverl\left(\A(x,\grl v)\right)-\abs {v}^{q-1}v
 = \diverl\left(\A(x,\grl u)\right)-\abs{u}^{q-1}u\qquad on\ \RN, $$
  and applying Theorem \ref{teo:comp} we conclude that $u\equiv v$.

  The remaining claim follows from Corollary \ref{cor:bd}.
\ep

\subsection{Some results for non \MPC\ operators }\label{sec:comp1}

Notice that the $p$-Laplacian operator with $p>2$ is not \MPC.
This fact it is easy to see by homogeneity consideration.

In this section we shall require that $p>2$ and  $\A:\Omega\times \R^l\to \R^l$
for all $\xi,\eta\in \R^l\setminus\{0\}$, $x\in \RN$ satisfies
\begin{equation}
  \label{mpcs}
  (\A(x,\xi)-\A(x,\eta))\cdot(\xi-\eta)\ge k_2
      \frac{\abs{\A(x,\xi)-\A(x,\eta)}^p}{(\abs \xi+\abs \eta)^{p(p-2)}}.
\end{equation}

\begin{example}
  Example of function $\A$ satisfying (\ref{mpcs})  is
$\A(x,\xi)=a(x)\abs \xi^{p-2}\xi$ with $a$ a bounded nonnegative 
function and $p\ge2$.
Indeed, if $\A(\xi)=\abs \xi^{p-2}\xi$, the following inequalities holds
\begin{eqnarray}
  \label{eq:p1}
  \abs{\xi^{p-2}\xi-\eta^{p-2}\eta}\le c_1 (\abs \xi+\abs \eta)^{p-1-\alpha}\abs{\xi-\eta}^\alpha\\
  (\xi^{p-2}\xi-\eta^{p-2}\eta )\cdot(\xi-\eta)\ge c_2 (\abs \xi+\abs \eta)^{p-\beta}\abs{\xi-\eta}^\beta\label{eq:p2}
\end{eqnarray}
with $\beta\ge \max\{p,2\}$ and $0\le \alpha\le\min\{1,p-1\}$.
See  \cite{Boccardo}.

Therefore choosing $\beta=p$ and $\alpha=1$ in (\ref{eq:p1}) and (\ref{eq:p2})
we have
\begin{eqnarray*}
   \frac{\abs{\A(x,\xi)-\A(x,\eta)}^p}{(\abs \xi+\abs \eta)^{p(p-2)}}=
  a(x)^p\frac{{\abs{\xi^{p-2}\xi-\eta^{p-2}\eta}}^p}{(\abs \xi+\abs \eta)^{p(p-2)}}
 \le a(x)^p c_1^p\abs{\xi-\eta}^p\\  
  \le a(x)^p \frac{c_1^p}{c_2} (\xi^{p-2}\xi-\eta^{p-2}\eta )\cdot(\xi-\eta)
     =  a(x)^{p-1} \frac{c_1^p}{c_2}(\A(x,\xi)-\A(x,\eta))\cdot(\xi-\eta).
\end{eqnarray*}
  Therefore, (\ref{mpcs}) is fulfilled with $k_2=\frac{c_1^p}{c_2}\norm a_\infty^{p-1}$.
\end{example}

We need of a version of Theorem \ref{teo:start} for operator
satisfying (\ref{mpcs}).
\begin{lemma}\label{teo:starts}  Let $\A:\Omega \times \R^l\to \R^l$ satisfy
  (\ref{mpcs}) with $p> 2$.
  Let $f,g\in L^1_{loc}(\Omega)$ and let $(u,v)$ be weak solution of
  (\ref{dis:comp0}). Set $w\decl(v-u)^+$ and let $s>0$. 
  If $(f-g) w\ge 0$ and 
  $w^{s+p'-1} \left(\abs{\grl u}+\abs{ \grl v}\right)^{p'(p-2)} \in L^1_{loc}(\Omega)$, then
  \be (f-g) w^s,\ \  (\A(x,\grl v)-\A(x,\grl u))\cdot \grl w\ w^{s-1}\in L^1_{loc}(\Omega),\label{regw+2}\ee
 and for any nonnegative $\phi\in\Cuno_0(\Omega)$ we have,
 \begin{eqnarray}  \int_\Omega (f-g) w^s \phi+c_1 
   \int_\Omega(\A(x,\grl v)-\A(x,\grl u))\cdot \grl w\ w^{s-1}\phi \\ \le
     c_2  \int_\Omega w^{s+p'-1}\left(\abs{\grl u}+\abs{ \grl v}\right)^{(p-2)p'} \frac{\abs{\grl\phi}^{p'}}{\phi^{p'-1}}, 
      \label{dis:ests2}\end{eqnarray}
  where $c_1=c_1(s,p,k_2),c_2(s,p,k_2)>0$ are suitable constants independent of $u$, $v$ and $\phi$.
\end{lemma}
\bp The proof is analogous to the proof of Theorem \ref{teo:start}.
  So we shall sketch it using the same notation.
  Applying Lemma \ref{th:katocompare} we have (\ref{tec1}) which,
  by using H\"older's inequality, (\ref{mpcs}) and Young's inequality, yields
\begin{eqnarray*}\lefteqn{\int_\Omega (f-g)\gamma(w) \phi+ \int_\Omega(\A_v-\A_u)\cdot\ \grl w\ \gamma'(w)\phi\le}\\
   &&\le\left(\int_\Omega\frac{\abs{\A_v-\A_u}^{p}}{\left(\abs{\grl v}+\abs{\grl u}\right)^{(p-2)p}}\gamma'(w)\phi\right)^{1/p}
    \left(\int_\Omega\frac{\gamma(w)^{p'}}{\gamma'(w)^{p'-1}}  \frac{\abs{\grl\phi}^{p'}}{\phi^{p'-1}}\left(\abs{\grl v}+\abs{\grl u}\right)^{(p-2)p'} \right)^{1/p'}\\
 && \le  \frac{\epsilon^{p}}{pk_2^p} 
\int_\Omega(\A_v-\A_u)\cdot \grl w\ \gamma'(w)\phi + 
   \frac{1}{p'\epsilon^{p'}}
  \int_\Omega\frac{\gamma(w)^{p'}}{\gamma'(w)^{p'-1}} \frac{\abs{\grl\phi}^{p'}}{\phi^{p'-1}}\,\left(\abs{\grl v}+\abs{\grl u}\right)^{(p-2)p'}.
  \end{eqnarray*}

  Next, constructing a sequence of $\gamma_n(t)$ approximating the function
  $t^s$ as made in the proof of  Theorem \ref{teo:start}, we conclude the proof.
\ep

\begin{theorem} Let $\A:\Omega \times\R^l\to \R^l$  $\A$ satisfy  (\ref{mpcs}) with $p>2$.
  Let $f,g\in L^1_{loc}(\Omega)$ and let $(u,v)$ be weak solution of
  (\ref{dis:comp0}). Set $w\decl(v-u)^+$ and let $s>0$. 
  If $(f-g) w\ge 0$ and 
  $w^{s(p-1)+1} \in L^1_{loc}(\Omega)$, then
  \be (f-g) w^s,\ \  (\A(x,\grl v)-\A(x,\grl u))\cdot \grl w\ w^{s-1}\in L^1_{loc}(\Omega)\label{regw1+}\ee
 and for any nonnegative $\phi\in\Cuno_0(\Omega)$ we have,
 \begin{eqnarray}
   \lefteqn{\int_\Omega (f-g) w^s \phi+c_1 
   \int_\Omega(\A(x,\grl v)-\A(x,\grl u))\cdot \grl w\ w^{s-1}\phi} \\ 
   && \quad \le
     c_2  \left(\int_A w^{s(p-1)+1} \frac{\abs{\grl\phi}^{p}}{\phi}\right)^{1/(p-1)}
     \left(\int_A\left(\abs{\grl u}+\abs{ \grl v}\right)^{p}\right)^{(p-2)/(p-1)} ,
      \label{dis:ests3}\end{eqnarray}
 where $A$ is the support of $\grl \phi$, and $c_1,c_2>0$ are 
 suitable constants independent of $u$, $v$ and $\phi$.
 
 In particular if there exist $c>0, q>0$ such that
 \begin{equation} \label{fmgmq}
   (f-g) (v-u)^+\ge c  \left((v-u)^+\right)^q,
 \end{equation}
  we have that $w^{q+s}\in L^1_{loc}(\Omega)$. 
  Moreover, $w^{q+1}\in L^1_{loc}(\Omega)$ and if $q>p-1$ we have
  \begin{equation}    \label{eq:intr}
    w^r\in L^1_{loc}(\Omega)\ \ for\  any\ \  r\in \left]0,\frac{q(p-1)-1}{p-2}\right[.  \end{equation}
\end{theorem}
\bp In order to apply Lemma \ref{teo:starts} we have to show that
   $w^{s+p'-1} \left(\abs{\grl u}+\abs{ \grl v}\right)^{p'(p-2)} \in L^1_{loc}(\Omega)$.

  Let  $\phi\in\Cuno_0(\Omega)$.
  An application of H\"older's inequality with exponent $z\decl \frac{p-1}{p-2}$
  implies
  \begin{equation}    \label{tec10}
   \int_\Omega w^{s+p'-1}  \frac{\abs{\grl\phi}^{p'}}{\phi^{p'-1}} \left(\abs{\grl u}+\abs{ \grl v}\right)^{(p-2)p'} \le \left(\int_Aw^{s(p-1)+1} \frac{\abs{\grl\phi}^{p}}{\phi}  \right)^{1/z'}
   \left(\int_A \left(\abs{\grl u}+\abs{ \grl v}\right)^{p}\right)^{1/z}.  
  \end{equation}
   Since $\abs{\gr u}, \abs {\gr v} \in L^p_{loc}(\Omega)$ and 
   $w^{s(p-1)+1} \in L^1_{loc}(\Omega)$ by hypotheses, we obtain the claim.
   Using (\ref{tec10}) in (\ref{dis:ests2}) we obtain (\ref{dis:ests3}).

  Now assuming (\ref{fmgmq}), we have that
   $w^{q+s}\in L^1_{loc}(\Omega)$.

  In order to complete the proof we begin observing that since 
  $w\in L^p_{loc}(\Omega)$, by choosing $s=1$, we have $u^{q+1}\in L^1_{loc}(\Omega)$.

  If $\frac{q(p-1)-1}{p-2}>q+1 $, that is if 
  $q>p-1$, we shall use a bootstrap argument.
  If $w^r\in  L^1_{loc}(\Omega)$, by the first part of the theorem,
  we have that $w^{h(r)}\in  L^1_{loc}(\Omega)$ with
  $h(r)\decl q+\frac{r-1}{p-1}$. Therefore, setting $r_0=q$ and 
  $r_{n+1}\decl h(r_n)$ we easily verify that the sequence $(r_n)_n$
  is increasing and it converges to $\frac{q(p-1)-1}{p-2}$. 
  This concludes the proof.
\ep

\begin{theorem}\label{teo:p2} Let $\A$ satisfy (\ref{mpcs}) with $p>2$.
  Let $(u,v)$ be a solution of (\ref{dis:cmp}) with    $q>1$.
  Assume that
 \be \left(\mint_{A_R} \abs {\grl v}^p \right)^{1/p}, \left(\mint_{A_R} \abs {\grl u}^p \right)^{1/p}\le c R^\theta\ \ for \ \ R\ large,\label{grp}\ee
with
\begin{equation}\label{grp2}
  \theta<\frac{1}{p-2}-\frac{Q}{p}
\end{equation}
  Then $v\le u$.
\end{theorem}
\bp Set $w\decl(v-u)^+$. From (\ref{wki}) we get that (\ref{fmgmq}) is satisfied, 
  and hence (\ref{eq:intr}) and (\ref{dis:ests3}) hold. 
  Therefore, fixing   $0<s<\frac{q-1}{p-2}$ we have that 
  $$\frac{q(p-1)-1}{p-2}> q+s>s(p-1)+1,$$
  and hence $w^{s(p-1)+1}\in L^1_{loc}(\RN)$.
  From (\ref{dis:ests3}) we have
  $$  \int_\rn w^{q+s} \phi \le
    c  \left(\int_A w^{s(p-1)+1} \frac{\abs{\grl\phi}^{p}}{\phi}\right)^{1/(p-1)}
    \left(\int_A\left(\abs{\grl u}+\abs{ \grl v}\right)^{p}\right)^{(p-2)/(p-1)}.$$
  By H\"older's inequality with exponent $x\decl\frac{q+s}{s(p-1)+1}$, and 
  choosing $\phi=\phi_R$ with $\phi_R$ as in the proof of Lemma \ref{lem:estR}, we have
  \begin{eqnarray}
    \left(  \int_{B_R} w^{q+s} \right)^{1-\frac1{x(p-1)}} &\le&
    \left(\int_A \frac{\abs{\grl\phi}^{px'}}{\phi^{x'}}\right)^{\frac1{x'(p-1)}}
    \left(\int_A\left(\abs{\grl u}+\abs{ \grl v}\right)^{p}\right)^{(p-2)/(p-1)}\nonumber\\
    &\le&    c R^{\frac{N-px'}{x'(p-1)}} R^{N\frac{p-2}{p-1}}\left(\mint_A\left(\abs{\grl u}+\abs{ \grl v}\right)^{p}\right)^{(p-2)/(p-1)}\nonumber\\
  &\le& c R^ t. \label{dis:tec13}
  \end{eqnarray}
  Here the last inequality follows form (\ref{eq:grh}), where 
  $$t\decl\frac{Q-px'}{x'(p-1)}+ Q\frac{p-2}{p-1} + \theta p\frac{p-2}{p-1}=
     \frac{Q}{x'(p-1)}+p \frac{p-2}{p-1} \left(\theta + \frac Qp-\frac{1}{p-2}\right).$$ 
  Since $$\lim_{s\to \frac{q-1}{p-2}^- }x'= 
  \lim_{s\to \frac{q-1}{p-2}^- }\frac{q+s}{q-1-s(p-2)}= +\infty,$$
  we can choose  $s$  so that $x'$ is large enough and $t<0$.
  By this choice, letting $R\to +\infty$ in (\ref{dis:tec13}), we obtain
  that $w\equiv 0.$ This proves our claim.
\ep

\begin{remark}
  Assumptions (\ref{grp}) and (\ref{grp2}) are obviously satisfied when
  looking for compact supported solutions.
\end{remark}
Further examples 
when the growth condition (\ref{grp}) holds are 
stated in the following.
\begin{proposition}\label{prop:reg} Let $q>p-1>0$ and $\A$ be \SPC. Let $h\in L^1_{loc}(\RN)$ 
  and let $u\in L^{q}_{loc}(\RN)\cap W^{1,p}_{loc}(\RN)$ be a weak solution of
  \begin{equation}    \label{eq:mainh}
     \diver (\A(x,\gr u))=\abs u^{q-1}u+h\quad on\ \ \RN.  \end{equation}
  Assume that  $u\, h\in L^1_{loc}(\RN)$. Then $u^{q+1}\in L^1_{loc}(\RN)$.
  Furthermore if $h\in L^{1+1/q}_{loc}(\RN)$, then for any $R>0$ we have
   $$ \int_{B_R} \abs u ^{q+1}  +  \int_{B_R} \abs{\grl u}^p\phi\le 
  c_4 R^{Q-p\frac{q+1}{q-p+1}} + c_5\int_{B_{2R}} \abs h^{\frac{q+1}{q}}  .$$
  In particular, if  there exists $\sigma \in \R$ such that 
  \begin{equation}    \label{eq:grh}
   \left(\mint_{B_R}\abs h^{1+1/q}\right)^{q/(q+1)} \le c R^\sigma\ \ for\ \ R\ large ,   \end{equation}
  then 
  \be \left(\mint_{B_R} \abs {\grl u}^p \right)^{1/p}\le c R^\theta\ \ for \ \ R\ large,\label{grp1}\ee
  where $$\theta\decl \max\left\{ \sigma \frac{q+1}{qp}, -\frac{q+1}{q+1-p}\right \}.$$
\end{proposition}
\bp We can use $u\phi$ as test function in (\ref{eq:mainh}).
  This follows by the fact that   $u\, h \in L^1_{loc}(\RN)$. 
  To see this  we argue as in Theorem \ref{teo:start} (see also \cite{dam-mit:kato}) obtaining
  $$ \int \abs u ^{q+1} \phi + \int \A(x,\gr u))\cdot \gr u\phi 
    \le \int \abs{\A(x,\gr u)} \abs{\gr \phi} \abs u +\int \abs h \abs u \phi.  $$
  By using the fact that $\A$ is \SPC, and by Young's inequality we get
  $$ \int \abs u ^{q+1} \phi + c_1 \int \abs{\gr u}^p\phi\le 
  c_2 \int \abs{\gr \phi}^p\abs{u}^p+\int \abs h \abs u \phi.  $$
  This implies that $u^{q+1}\in L^1_{loc}(\RN)$.

  Next, assume that $h\in L^{1+1/q}_{loc}(\RN)$.
  By Young's inequality with exponents $x\decl\frac{q+1}{p}$ and
  $y\decl {q+1}$, it follows that
   $$ c_3\int \abs u ^{q+1} \phi + c_1 \int \abs{\gr u}^p\phi\le 
  c_4\int \abs{\gr \phi}^{px'} + c_5\int \abs h^{y'} \phi.  $$
  Next by choosing $\phi=\phi_R$ as in Lemma \ref{lem:estR}, 
  from the assumption  on $h$ we get
  $$\mint_{B_R} \abs{\gr u}^p\le c (R^{-px'}+ R^{\sigma\frac{q+1}{q}})
  \le c R^{\max \{-px', \sigma\frac{q+1}{q} \}}. $$
\ep

\begin{remark}
  The assumption $u\, h\in L^1_{loc}(\RN)$ is obviously satisfied for instance if 
  $u\in L^{q+1}_{loc}(\RN)$ and  $h\in L^{1+1/q}_{loc}(\RN)$.
\end{remark}

\begin{corollary} Let $q>p-1>1$. Let $\A$ be \SPC\ satisfying (\ref{mpcs}) and
  let $h\in L^{1+1/q}_{loc}(\RN)$ be 
  such that  (\ref{eq:grh}) holds for $\sigma \in \R$ and  
  $$ \max\left\{ \sigma \frac{q+1}{qp}, -\frac{q+1}{q+1-p}\right \}<\frac{1}{p-2}-\frac{N}{p}. $$
  Then problem (\ref{eq:mainh}) has at most one solution on the class
  $  L^{q+1}_{loc}(\RN)\cap W^{1,p}_{loc}(\RN)$.
\end{corollary}
\begin{remark}
The possible weak solutions belong to the space  $ L^{q+1}_{loc}(\RN)$
in the following cases.
\begin{enumerate}
\item Since $u\in W^{1,p}_{loc}(\RN)$ then, by Sobolev's embedding, 
  $u\in L^{q+1}_{loc}(\RN)$ provided $q\le \frac{N(p-1)+p}{N-p}$.
\item If $p\ge N$.
\item If $h\in L^\infty_{loc}(\RN)$ then $u\, h\in L^1_{loc}(\RN)$. By  
  Proposition \ref{prop:reg} it follows that $ u\in L^{q+1}_{loc}(\RN)$
\item If $h\in L^\infty(\RN)$ then from Theorem \ref{cor:bd} it follows that $u\in  L^\infty(\RN)$.
  In particular Corollary \ref{cor:lim} holds.
\end{enumerate}
\end{remark}
\begin{corollary} \label{cor:lim} Let $q>p-1>1$. 
  Let $\A$ be \SPC\ satisfying (\ref{mpcs}) and  $h\in L^{\infty}(\RN)$.
  If 
$$ 2<p<\frac{2Q}{Q-1},$$
  then problem (\ref{eq:mainh})  has at most one weak solutions.
\end{corollary}

\begin{corollary} Let $h\in L^1_{loc}(\RN)$, $q>p-1$ and $p>2$.
  Then the problem 
  $$-\Delta_{L,p}u+ \abs u ^{q-1}u =h,$$
  has at most one weak solution $u$
  satisfying, $\displaystyle{\int_{\rn} \abs{\grl u}^p<+\infty}$.
\end{corollary}
\bp
  It is enough to choose $\theta \decl -Q/p$ in 
  Theorem \ref{teo:p2}.
\ep

Now requiring stronger assumptions on the behavior of the gradient of the 
solutions, we have the following.

\begin{theorem}\label{teo:p2g} Assume that $\A$ satisfies condition (\ref{mpcs}) with $p>2$.
  Let $(u,v)$ be a solution of (\ref{dis:cmp}) with  $q>1$.

  Let $\theta <\frac{1}{p-2}$ and assume
  that there exists $\alpha> \frac{N(p-2)}{1-\theta(p-2)}$ such that
  \begin{eqnarray}
   \alpha \frac{q(p-1)-1}{p-2}>qp,\label{57}\\
  \left((v-u)^+\right)^{\alpha \frac{q(p-1)-1 }{p(p-2)}}\in L^1(A_R) \ \  for \ \ R\ \  large,\\
   \left(\mint_{A_R} \abs {\grl v}^\alpha \right)^{1/\alpha}, \left(\mint_{A_R} \abs {\grl u}^\alpha \right)^{1/\alpha}\le c R^\theta.\label{grcr}
  \end{eqnarray}

  Then $v\le u$ a.e. on $\RN$.
\end{theorem}
\bp Let $w\decl(v-u)^+$. By (\ref{wki}) and
    Lemma \ref{teo:starts}, 
   applying H\"older's inequality to (\ref{dis:ests2}) with exponent
   $x\decl \frac{q+s}{s-p'+1}$ where $s>0$, we have
   \be \int_\Omega w^{q+s}\phi \le c_2 \left( \int_\Omega w^{q+s}\phi\right)^{1/x}
    \left(\int_\Omega \left(\abs{\grl u}+\abs{ \grl v}\right)^{(p-2)p'x'} \frac{\abs{\grl\phi}^{p'x'}}{\phi^{p'x'-1}}\right)^{1/x'}.
\label{tec4} \ee

  By choosing $\phi=\phi_R$  as in the proof of Lemma \ref{lem:estR},
  it follows that
  \be \int_{B_R} w^{q+s} \le c_2 c(\phi_1) R^{Q-p'x'}
\mint_{A_R} \left(\abs{\grl u}+\abs{ \grl v}\right)^{(p-2)p'x'}. \label{tec5}\ee
  Next, we observe that we must choose $s>0$ such that  all the integrals in (\ref{tec4}) are well defined.
  Indeed, by choosing $s\decl \frac{\alpha}{p'(p-2)}(q-p'+1)-q$, from (\ref{57}) it follows that $s>0$.
  We observe that since  $(p-2)p'x'= \alpha$ and  
  $w^{q+s}= w^{\frac{\alpha}{p'(p-2)}(q-p'+1)}\in L^1_{loc}(\RN)$,
  the integrals in (\ref{tec5}) and  (\ref{tec4}) are well defined.
  Next, from (\ref{tec5}) and (\ref{grcr}), we obtain
  \be \int_{B_R} w^{q+s} \le c R^\gamma \qquad where\ \ \gamma={Q-p'x' +\theta (p-2)p'x'}. \label{tec6}\ee

  Finally, we observe that since
$$\gamma =Q- \alpha \frac{1}{p-2}+\theta \alpha<Q-(\frac{1}{p-2}-\theta)
  \frac{Q(p-2)}{1-\theta(p-2)}=0,$$
  by letting $R\to+\infty$ in (\ref{tec6}), the claim follows.
\ep

As a final remark we note that the knowledge of
a pointwise estimate on the gradient of the solutions on an exterior domain
of the type
$$ \abs{\gr u(x)} \le c \abs x^\theta\quad \ for \ \ \abs x\ \ large,$$
where $\theta<1/(p-2)$, then Theorem \ref{teo:p2g} applies and yields
Theorem \ref{14}. 

\noindent{\bf Proof of Theorem \ref{14}.} The growth assumption on $\gr u$ implies that
 (\ref{grcr}) holds for any $\alpha>0$. By the Sobolev embedding theorem
$u \in L^r(A_R)$ for any $r>0$ and for $R$ large. 
Therefore by choosing $\alpha$ large enough the claim follows from 
Theorem \ref{teo:p2g}.
\ep

\appendix
\section{Inequalities and  \MPC\ Operators}\label{appMPC}

Here, we shall prove some fundamental elementary inequalities 
that we use throughout the paper.
Very likely these inequalities are well known, nevertheless for 
completeness we shall include their proof here.

In what follows we shall assume that $\A$ has the form
$$\A(x,\xi)=A(\abs\xi)\xi,$$
where $\A:\R_+\to \R$. We set $\phi(t)\decl A(t)t$.

\begin{theorem}\label{teo:curvmpc} 
 Let $A$ be nonincreasing and bounded function such that
 \be \phi(0)=0,\ \ \phi(t)>0 \mathrm{\ for\ } t>0, \phi\textrm{\ is\ nondecreasing}.
  \label{hyp:cr} \ee
  Then $\A$ is \MPC\ with $p=2$.    
\end{theorem}

\begin{theorem}\label{teo:plapmpc} Let $1<p\le 2$. 
  Let $\phi$ be increasing, concave function satisfying
  (\ref{hyp:cr}) and such that there exist positive constants $c_p,c_\phi>0$ such that
 \be \phi(t)\le c_p t^{p-1} 
     \label{hyp:fip} \ee  and
\be \phi'(s) s \le c_\phi \phi(s).\label{hyp:crp}\ee
  Then $\A$ is \MPC.
\end{theorem}

\begin{remark} We notice that \ref{hyp:fip} is necessary condition
  for $\A$ to be an \MPC\ operator.
 Indeed, if $\A$ is \MPC, by taking $\eta =0$, then it follows that
   $\A$ is \WPC, and  (\ref{hyp:fip}) holds by H\"older inequality.
\end{remark}
We set
$$I\decl (\A(\xi)-\A(\eta))\cdot (\xi-\eta), \qquad
    J \decl \abs{\A(\xi)-\A(\eta)}.$$ 

Our goal is to prove that there exists a constant $c>0$ 
such that $I^{p-1}\ge c J^p $.

We set $t\decl \abs \xi$, $s\decl \abs \eta$ and let $\theta$ be such that
$\xi\cdot\eta=\theta \abs \xi  \abs \eta = \theta t s$. 
Hence $\theta\in [-1.1]$, $t,s >0$. Moreover by symmetry we can assume
that $s\ge t$.

We rewrite $I$ and $J$ as
\begin{eqnarray*}
   I &=& A(\abs \xi) \abs \xi ^2 + A(\abs \eta) \abs \eta ^2
    -A(\abs \xi)\xi\cdot\eta -A(\abs \eta)\xi\cdot\eta\\
    &=&
    \phi(t) t +\phi(s) s -\phi(t) s \theta - \phi(s) t\theta, \\
 J^2 &=& \phi^2(t)+\phi^2(s)-2\phi(t)\phi(s)\theta .
\end{eqnarray*}

\begin{remark}
  From (\ref{hyp:cr}) we deduce that: if $I=0$ then $\phi(t)=\phi(s)$, 
 $\theta =1$, and $J=0$. Indeed, assuming $s \ge t$
$$I = \phi(t) (t-s\theta) +\phi(s) (s-t\theta)\ge
  \phi(t) (t-s\theta) + \phi(t) (s-t\theta)=\phi(t)(1-\theta)(t+s)\ge0 $$
Therefore, if $I=0$, then $\theta =1$ or $\phi(t)=0$. If $\phi(t)=0$ then
  $t=0$ and hence (since $I=0$) also $s=0$.
  If $\theta=1$, then we have
  $ 0= I=  (\phi(s)-\phi(t)) (s-t)$ and hence the claim follows.

  We notice that if $\phi$ is increasing, then $I=0$ implies also that $t=s$.

  Therefore, in order to prove that $\A$ is \MPC, we restrict
  ourselves to the case $s>t>0$.
\end{remark}

\noindent{\bf Proof of Theorem \ref{teo:curvmpc}.} 
  Set
$$ I_1:= \frac{I}{\phi(t)t}=1+\frac{\phi(s)s}{\phi(t)t}-\theta\frac{\phi(s)}{\phi(t)}
  -\theta \frac{s}{t}
$$
$$ J_1:= \frac{J^2}{\phi^2(t)}=1+\frac{\phi^2(s)}{\phi^2(t)}-2\theta\frac{\phi(s)}{\phi(t)}
$$ 
We have
\begin{eqnarray*}
   I_1-J_1&=&\frac{\phi(s)s}{\phi(t)t}-\frac{\phi^2(s)}{\phi^2(t)}+
  \theta\frac{\phi(s)}{\phi(t)} -\theta \frac{s}{t}=
  \frac{\phi(s)}{\phi(t)}\left(\frac st-\frac{\phi(s)}{\phi(t)}\right)
  -\theta\left(\frac st-\frac{\phi(s)}{\phi(t)}\right)\\
  &=&\left(\frac st-\frac{\phi(s)}{\phi(t)}\right)\left(\frac{\phi(s)}{\phi(t)}-\theta \right)=\frac st\left(1-\frac{A(s)}{A(t)}\right)\left(\frac{\phi(s)}{\phi(t)}-\theta \right).
\end{eqnarray*}

Since $s> t>0$,  $\phi$ is nondecreasing and 
$A$ is nonincreasing it follows that  $I_1-J_1\ge 0$.

Therefore
$$ I=\phi(t) t I_1\ge \phi(t) t J_1 =\frac{t}{\phi(t)} J^2=\frac{1}{A(t)} J^2\ge
  \frac{1}{\norm A_\infty} J^2=cJ^2,$$
that is the claim. \ep

\noindent{\bf Proof of Theorem \ref{teo:plapmpc}.}
  Our goal is to show that  $ I^{p-1} / J^p\ge const>0$.

 We have
  \begin{eqnarray*}
    \frac{I^{p-1}}{ J^p}=\frac{s^{p-1}}{\phi(s)}\ \frac{(\frac{I}{\phi(s)s})^{p-1}}{\left(\frac{J^2}{\phi^2(s)}\right)^{p/2}}=\frac{s^{p-1}}{\phi(s)}\ 
     \frac{\left(1+\frac{\phi(t)t}{\phi(s)s}-\theta \frac{\phi(t)}{\phi(s)}-\theta\frac ts\right)^{p-1}}
     {\left(1+ \frac{\phi^2(t)}{\phi^2(s)}-2 \theta \frac{\phi(t)}{\phi(s)}\right)^{p/2}}\ge c F(t,s,\theta),
  \end{eqnarray*}
where
$$ F(t,s,\theta)\decl \frac{\left(1+\frac{\phi(t)t}{\phi(s)s}-\theta \frac{\phi(t)}{\phi(s)}-\theta\frac ts\right)^{p-1}}
     {\left(1+ \frac{\phi^2(t)}{\phi^2(s)}-2 \theta \frac{\phi(t)}{\phi(s)}\right)^{p/2}}.
$$
In order to prove the claim it is enough to show that $F$ is uniformly positive
  for $\theta\in [-1,1]$ and $s>t>0$.
  Since $\phi$ is nondecreasing, setting 
  $$\alpha\decl \frac{\phi(t)}{\phi(s)},\quad z\decl \frac ts,$$
  it is enough to prove that
  $$G(\alpha,z,\theta)\decl \frac{\left( 1+\alpha z-\alpha\theta-z \theta \right)^{p-1}}{\left(1+\alpha^2-2\theta\alpha\right)^{p/2}}
  $$
 is uniformly positive for 
 $(\alpha,z,\theta)\in D\decl [0,1]\times [0,1]\times [-1,1]\setminus\{1,1,1\}$.
 The function $G$ is well defined in $D$. Indeed the denominator vanishes
 if and only if  $J=0$ that is $\phi(t)=\phi(s)$ and $\theta =1$ that is
  $\alpha=z=\theta=1$.
  On the other hand the numerator of $G$ vanishes if $I=0$, that is
  if $\alpha=z=\theta=1$. Therefore, $G$ is strictly positive on $D$.

  Moreover taking into account that $\phi$ is concave we have
 $$\frac{\phi(s)-\phi(t)}{s-t}\le
   \phi'(s), $$
  which, together with (\ref{hyp:crp}), yields
 $$\frac{\phi(s)-\phi(t)}{s-t}\ \frac{s}{\phi(s)}\le 
   \phi'(s)\frac{s}{\phi(s)} \le c_\phi. $$
  That is 
  \be 1-\alpha\le c_\phi (1-z).\label{dis:constr}\ee

  Now the claim will follows by proving that
  \be \liminf_{\alpha,z,\theta\to 1} G(\alpha,z,\theta)>0.\label{limg}\ee
  Here the $\liminf$ is computed for $(\alpha,z,\theta)\in D$ and under the constrained  (\ref{dis:constr}).
  Introducing
  $$H(a,b,e)\decl \frac{\left(ab +2e -ae -be \right)^{p-1}}
    {\left(a^2+2e-2ae \right)^{p/2}}$$
  and setting
  $a\decl 1-\alpha$, $b\decl 1-z$, $e\decl 1-\theta$, 
  (\ref{limg}) is equivalent to
  $$ \liminf_{a,b,e \to 0^+} H(a,b,e)>0$$
  with $a\le c_\phi b$.

  We shall argue by contradiction. Let $a_n,b_n,e_n$ be three infinitesimal sequences 
  such that $H(a_n,b_n,e_n)\to 0$.
  Since $a_ne_n$ and $b_ne_n$   are  infinitesimal sequences of order greater then $e_n$, we have that
  $$ 0 = \lim_nH(a_n,b_n,e_n) = \lim_n \frac{\left(a_nb_n +2e_n \right)^{p-1}}
    {\left(a_n^2+2e_n \right)^{p/2}}.$$
  Taking into account that $ a_n\le c_\phi b_n$,
  we have
  $$ 0=\lim_n \frac{\left(a_nb_n +2e_n \right)^{p-1}}
    {\left(a_n^2+2e_n \right)^{p/2}}\ge
    \liminf_n \frac{\left(a_n^2/c_\phi  +2e_n \right)^{p-1}}
    {\left(a_n^2+2e_n \right)^{p/2}} >0, $$
  because $a_n^2/c_\phi  +2e_n $ and $a_n^2+2e_n $ are infinitesimal of the same order
  and $p\le 2$. This contradiction concludes the proof.
\ep

\begin{remark} If $\A$ has the form
  $$\A(x,\xi)=a(x)A(\abs\xi)\xi,$$
  then the above Theorems \ref{teo:curvmpc} and \ref{teo:plapmpc} hold
  provided $a\in L^\infty(\RN)$ and it is positive a.e..
\end{remark}

\section{Carnot Groups}\label{appCarnot}

We quote some  facts on Carnot groups  and refer the interested reader to  
\cite{bon-lan-ugu07} 
  for  more detailed information on this subject.

A Carnot group is a connected, simply connected, nilpotent Lie
group $\G$ of dimension $N$ with graded Lie algebra ${\cal
G}=V_1\oplus \dots \oplus V_r$ such that $[V_1,V_i]=V_{i+1}$ for
$i=1\dots r-1$ and $[V_1,V_r]=0$. Such an integer $r$ is called the
\emph{step} of the group.
 We set $l=n_1=\dim V_1$, $n_2=\dim V_2,\dots,n_r=\dim V_r$.
A  Carnot group $\G$ of dimension $N$ can be identified, up to an
isomorphism, with the structure of a \emph{homogeneous Carnot
Group} $(\RN,\circ,\delta_R)$ defined as follows; we
identify $\G$ with $\RN$ endowed with a Lie group law $\circ$. We
consider $\RN$ split in $r$ subspaces
$\RN=\RR^{n_1}\times\RR^{n_2}\times\cdots\times\RR^{n_r}$ with
$n_1+n_2+\cdots+n_r=N$ and $\xi=(\xi^{(1)},\dots,\xi^{(r)})$ with
$\xi^{(i)}\in\RR^{n_i}$. 
We shall assume that for any $R>0$ the dilation
$\delta_R(\xi)=(R\xi^{(1)},R^2 \xi^{(2)},\dots,R^r \xi^{(r)})$ 
is a Lie group automorphism.
The Lie algebra of
left-invariant vector fields on $(\RN,\circ)$ is $\cal G$. For
$i=1,\dots,n_1=l $ let $X_i$ be the unique vector field in $\cal
G$ that coincides with $\partial/\partial\xi^{(1)}_i$ at the
origin. We require that the Lie algebra generated by
$X_1,\dots,X_{l}$ is the whole $\cal G$.

We denote with $\grl$ the vector field $\grl\decl(X_1,\dots,X_l)^T$
and we call it \emph{horizontal vector field} and by  $\diverl$ the formal 
adjoint on $\grl$, that is (\ref{eq:defdiv}).
Moreover, the vector
fields $X_1,\dots,X_{l}$ are homogeneous of degree 1 with respect
to $\delta_R$ and in this case
$Q= \sum_{i=1}^r i\,n_i=  \sum_{i=1}^r i\,\mathrm{dim}V_i$ is called the 
\emph{homogeneous dimension} of $\G$.
The \emph{canonical sub-Laplacian} on $\G$ is the
second order differential operator defined by
$$\Delta_G=\sum_{i=1}^{l} X_i^2=\diverl(\grl\cdot )$$ and for $p>1$ the
$p$-sub-Laplacian operator is
$$\Delta_{G,p} u\decl \sum_{i=1}^{l} X_i(\abs{\grl u}^{p-2}X_iu)=
   \plap{u}.$$
Since $X_1,\dots,X_{l}$
generate the whole $\cal G$, the sub-Laplacian $\Delta_G$ satisfies the
H\"ormander hypoellipticity  condition. 

In this paper $\nabla$ and $|\cdot|$ stand  respectively for the
usual gradient in $\RN$
and the Euclidean norm.

Let $\mu\in \C(\RR^N;\RR^l)$ be a matrix 
$\mu\decl(\mu_{ij})$, $i=1,\dots,l$, $j=1,\dots,N$.
For $i=1,\dots,l$,  let $X_i$ and its formal adjoint $X_i^*$
be defined as
\begin{equation} X_i\decl\sum_{j=1}^N \mu_{ij}(\xi)\frac{\partial}{\partial \xi_j},
\qquad X_i^*\decl-\sum_{j=1}^N \frac{\partial}{\partial\xi_j}\left(\mu_{ij}(\xi)\cdot\right),
   \label{mu2}\end{equation}
and let $\grl$ be the vector field defined by
$\grl\decl (X_1,\dots,X_l)^T=\mu\nabla$
and $\grl^*\decl(X_1^*,\dots,X_l^*)^T$.

For any vector field $h=(h_1,\dots,h_l)^T\in\Cuno(\Omega,\RR^l)$, we shall use
the following notation $ \diverl(h)\decl\diver(\mu^T h)$,
that is
\begin{equation}  \label{eq:defdiv}
   \diverl(h)=-\sum_{i=1}^l X_i^*h_i=-\grl^*\cdot h.\end{equation}

An assumption that we shall made (which actually is an assumption on the matrix
$ \mu$) is that the operator
$$  \Delta_Gu =\diverl(\grl u)$$
is a canonical sub-Laplacian on a Carnot group 
(see below for a more precise meaning).
The reader, which is not acquainted with these structures, 
can think to the special case 
of $\mu =I$, the identity matrix in $\RN$, that is the usual Laplace operator
in Euclidean setting.

A nonnegative continuous function $S:\RN\to\RR_+$ is called a 
\emph{homogeneous norm} on { $\G$}, if 
$S(\xi^{-1})=S(\xi)$, $S(\xi)=0$ if and only if $\xi=0$, and it is
homogeneous of degree 1 with respect to $\delta_R$ (i.e.
$S(\delta_R(\xi))=R S(\xi)$).
A homogeneous norm $S$ defines on $\G$ a \emph{pseudo-distance} defined as
$d(\xi,\eta)\decl S(\xi^{-1}\eta)$, which
in general is not a distance.
If $S$ and $\tilde S$ are two homogeneous norms, then they are equivalent,
that is, there exists a constant
$C>0$ such that $C^{-1}S(\xi)\le \tilde S(\xi)\le CS(\xi)$.
Let $S$ be a homogeneous norm, then there exists a constant
$C>0$ such that $C^{-1}\abs\xi\le S(\xi)\le C\abs\xi^{1/r}$,
for $S(\xi)\le1$.
An example of homogeneous norm is
$  S(\xi)\decl\left(\sum_{i=1}^r\abs{\xi_i}^{2r!/i}\right)^{1/2r!}.$

Notice that if $S$ is a homogeneous norm differentiable a.e., 
then $\abs{\grl S}$ is homogeneous of degree 0 with respect to 
$\delta_R$; hence $\abs{\grl S}$ is bounded.

We notice that in a Carnot group, the Haar measure coincides with the Lebesgue measure.

\medskip

Special examples of Carnot groups are the
  Euclidean spaces $\R^Q$.
  Moreover, if $Q\le 3$ then any Carnot group is the ordinary Euclidean
  space $\RR^Q$.

 The simplest nontrivial example of a Carnot group
  is the Heisenberg group $\hei^1=\RR^3$.
  For an integer $n\ge1$, the Heisenberg group $\hei^n$ is defined as follows:
  let $\xi=(\xi^{(1)},\xi^{(2)})$ with
  $\xi^{(1)}\decl(x_1,\dots,x_n,y_1,\dots,y_n)$ and $\xi^{(2)}\decl t$.
  We endow $\RR^{2n+1}$ with  the group law
$\hat\xi\circ\tilde\xi\decl(\hat x+\tilde x,\hat y+\tilde y,\hat t+ \tilde t+2\sum_{i=1}^n(\tilde x_i\hat y_i-\hat x_i \tilde y_i)).$
We consider the vector fields
\[X_i\decl\frac{\partial}{\partial x_i}+2y_i\frac{\partial}{\partial t},\
        Y_i\decl\frac{\partial}{\partial y_i}-2x_i\frac{\partial}{\partial t},
        \qquad\mathrm{for\ } i=1,\dots,n, \]
and the associated  Heisenberg gradient
$ \grh\decl (X_1,\dots,X_n,Y_1,\dots,Y_n)^T$.
The Kohn Laplacian $\lh$ is then the operator defined by
$\lh\decl\sum_{i=1}^nX_i^2+Y_i^2.$
The family of dilations is given by
$\delta_R(\xi)\decl (R x,R y,R^2 t)$ with homogeneous dimension
$Q=2n+2$.
In ${\hei^n}$ a canonical homogeneous norm is defined as
$\abs{\xi}_H\decl \left(\left(\sum_{i=1}^n x_i^2+y_i^2\right)^2+t^2\right)^{1/4}.$

\subsection*{Acknowledgments}
\medskip
{This work is  supported by the Italian MIUR National Research Project: Quasilinear Elliptic Problems and Related Questions.}

\bibliographystyle{amsplain}

\end{document}